\documentclass{article}
\usepackage[utf8]{inputenc}
\usepackage{amsmath}
\usepackage{booktabs}
\usepackage{arydshln}
\usepackage{mdframed}
\usepackage{lipsum}
\usepackage{graphicx}

\newtheorem{theorem}{Theorem}
\newtheorem{conjecture}{Conjecture}
\newtheorem{definition}{Definition}
\newtheorem{claim}{Claim}
\newtheorem{question}{Question}
\newtheorem{remark}{Remark}

\title{\textbf{How Low Can You Go?} \\
New Bounds on the Biplanar Crossing Number \\of Low-dimensional Hypercubes}
\author{Gregory J. Clark\thanks{University of South Carolina. Columbia, SC, USA.}, Gwen Spencer\thanks{Smith College.
Northampton, MA, USA. \textit{
This material is based upon work supported by the National Science Foundation under Grant Number DMS 1641020 while both authors attended the ``Beyond Planarity: Crossing Numbers of Graphs" workshop (organized by the American Mathematical Society), and
under Grant No. DMS-1440140 while G. Spencer was in residence at the
Mathematical Sciences Research Institute in Berkeley, California, during the Fall 2017 semester.}}}

\date{June 2017}

\begin{document}

\maketitle

\begin{abstract}
In this note we provide an improved upper bound on the biplanar crossing number of the 8-dimensional hypercube. The $k$-planar crossing number of a graph $cr_k(G)$ is the 
number of crossings required when every edge of $G$ must be drawn in  
one of $k$ distinct planes.
It was shown in \cite{Cza} that $cr_2(Q_8) \leq 256$ which we improve to $cr_2(Q_8) \leq 128$.  Our approach highlights the relationship between symmetric drawings and the study of $k$-planar crossing numbers.  We conclude 
with several open questions concerning this relationship.
\end{abstract}

\section{Introduction}

The traditional \textit{crossing number} of a graph $G=(V,E)$, denoted by $cr(G)$, is the minimum number of edge crossings required to draw $G$ in the 2-dimensional Euclidean plane. To study printed circuit boards, Owens \cite{Owe} generalized the question: what is the minimum number of edge crossings required by a drawing that is allowed to carefully divide the edges of $G$ among two different 2-dimensional Euclidean planes?  Since then the definition has been extended to $k \geq 2$ planes \cite{Cza}.

Suppose that $E$ is partitioned into $k$ disjoint subsets, $E_1,E_2,...,E_k$, and let $G_i=(V, E_i)$. Each $G_i$ has some crossing number $cr(G_i)$. Suppose further that $G_i$ will be drawn in the $i$th plane from a set of $k$ distinct planes. The \emph{$k-planar$ crossing number of $G$}, denoted $cr_k(G)$ is then the minimum of
\[cr(G_1)+cr(G_2)+...+cr(G_k)\]
over all partitions of the edge set $E$. 

Trivially, letting $E_1=E$ shows that $cr_k(G)\leq cr(G)$. The question remains: given the freedom to consider any partition of $G$'s edges among $k$ disjoint planes, how low can we drive the number of required crossings?

A significant challenge in designing a crossing-minimizing $k$-planar drawing of $G$ is that, even for quite simple $G_i$, $cr(G_i)$ could be unknown. For example: for $Q_4$, the 4-dimensional hypercube, it is known that $cr(Q_4)=8$; however, the exact value of $cr(Q_d)$ is unknown for $d>4$ \cite{Far}.

The previous upper bound $cr_2(Q_8) \leq 256$ was given by a construction of Czabarka, S\'ykora, Sz\'ekely, and Vr\'to in \cite{Cza}. Czabarka et al. give a general construction for an upper bound on $cr_2(Q_d)$ that achieves 256 crossings when $d=8$. Their approach specifies a bi-planar partition of the edges of $Q_8$  based on a set of lower-dimensional hypercube subgraphs. Their upper bound is minimized when these hypercube subgraphs are as-uniform-as-possible in size. In particular, for $Q_8$ their construction specifies sixteen disjoint $Q_4$ subgraphs in Plane 1 and a further sixteen disjoint $Q_4$ subgraphs in Plane 2.  Recall that $cr(Q_4)=8$, so drawing each disjoint copy of $Q_4$ optimally yields \[cr_2(Q_8) \leq 16\times2\times8=256.\]

We now present our main result which improves on the the best known upper bound of $cr(Q_8)$ by a factor of 2.
\begin{theorem}\label{Tmain}
There exists a 2-planar drawing of the 8-dimensional hypercube with 128 crossings so that $cr_2(Q_8)\leq 128$.
\end{theorem}

%We will demonstrate several different biplanar drawings of $Q_8$ that achieve 128 crossings. Our first drawing is explicitly described in this section, and two additional drawings appear in the appendix.
\section{A biplanar drawing of $Q_8$ with 128 crossings}

To prove Theorem \ref{Tmain}, we  provide a biplanar drawing of $Q_8$ with 128 crossings.  We improve the previous construction by plane-swapping edges to give a net reduction in total edge crossings.  Our drawing consists of graphs $G_1$ and $G_2$ in Plane 1 and 2 respectively such that $G_1 \cong G_2$ where $cr(G_i) \leq 64$.  We found several distinct bi-planar drawings of $Q_8$ with exactly 128 crossings which satisfy these conditions. For ease of exposition, we present a highly symmetric drawing.

We define a \emph{depleted $n$-dimensional hypercube} to be a graph whose vertex set is $V(Q_n)$ and will refer to such graphs as \textit{depleted $n$-cubes}.  We will make use of depleted 5-cubes.  To this end we introduce the following partition $V(Q_4) := C_1 \sqcup C_2$ where 
\begin{align*}
C_1&:=\{0000,1000,0010,1010,0011,1011,0001,1001\} \\
 C_2 &:= \{0111,1111,0101,1101,0100,1100,0110,1110\}.
 \end{align*}
 
 For ease of notation, we denote $\hat c \in C_1$ and $\check c \in C_2$.  Moreover, we let $b \in \{0,1\}$ represent the usual binary-bit.  Maintaining the notation of \cite{Cza} we refer to each node of $Q_8$ by a length-8 binary string from $\{0,1\}^8$.  Given two binary strings $s_1$ and $s_2$ we write $s_1s_2$, or $s_1-s_2$ for readability, to be the usual string concatenation. 
 
In our construction, each plane contains 512 edges, and furthermore, $G_1$ and $G_2$ are isomorphic.  For exposition, suppose that we initially have a Plane 0 which contains all the edges and vertices of $Q_8$.  Further suppose that there exist Planes 1 and 2 which each initially contain the vertices of $Q_8$ and no edges.  We move every edge from Plane 0 to either Plane 1 or Plane 2 to create our biplanar partition.  In the following table, we describe explicitly the 512 edges we add to Plane 1.

Consider the set of pairs \[P_1:=\{(0000,1000), (0010, 1010), (0011, 1011), (0001, 1001)\} \subset \binom{C_1}{2}.\]  For $(\hat c_1, \hat c_2) \in P_1$ define the \emph{depleted 5-cube of Type 1}, denoted $D_1(\hat c_1, \hat c_2)$, according to Table \ref{Type1}.

\begin{table}[h!]
\centering
\begin{tabular}{|c|c|c|c|}
    \hline
    \multicolumn{4}{|c|}{$E(D_1(\hat c_1, \hat c_2))$ for $\hat c \in (\hat c_1, \hat c_2) \in P_1$}\\
    \hline
   %outer square in drawing type 1:
  $(\hat{c}-b000, \hat{c}-b001)$& $(\hat{c}-b000, \hat{c}-b100)$& $ (\hat{c}-b100, \hat{c}-b101)$&$(\hat{c}-b001, \hat{c}-b101)$\\ \hline
   %inner square in drawing type 1
   $(\hat{c}-b010, \hat{c}-b011) $&$ (\hat{c}-b010, \hat{c}-b110)  $& $(\hat{c}-b110, \hat{c}-b111) $& $(\hat{c}-b011, \hat{c}-b111)$\\ \hline
   % diagonals in drawing type 1
   $(\hat{c}-b000, \hat{c}-b010)$ & $(\hat{c}-b001, \hat{c}-b011)$  &$(\hat{c}-b100, \hat{c}-b110)$ &$(\hat{c}-b101, \hat{c}-b111)$\\ \hline
   %Q_3 to Q_3 edges
       $(\hat{c}-0101, \hat{c}-1101)$ &$(\hat{c}-0111, \hat{c}-1111)$ &$(\hat{c}-0110, \hat{c}-1110)$& $(\hat{c}-0100, \hat{c}-1100)$\\
       \hline
%    $(0000-\hat c, 1000-\hat c)$ & $(0010-\hat c, 1010-\hat c)$ & $(0011-\hat c, 1011-\hat c)$ & $(0001 -\hat c , 1001-\hat c)$\\
  $(\hat c_1-0000, \hat c_2-0000)$ & $(\hat c_1-0100, \hat c_2-0100)$ & $(\hat c_1-1100, \hat c_2-1100)$ & $(\hat c_1-1000, \hat c_2-1000)$\\ \hline
    $(\hat c_1-1001, \hat c_2-1001)$ & $(\hat c_1-1101, \hat c_2-1101)$ & $(\hat c_1-0101, \hat c_2-0101)$ & $(\hat c_1-0001, \hat c_2-0001)$\\
  \hline
\end{tabular}
\caption{Table of the 64 edges of \textit{depleted 5-cubes of Type 1}.}
\label{Type1}
\end{table}

\begin{figure}[h!]
\label{F:Type1}
\hspace*{-17mm}\includegraphics[trim = 0mm 0mm 0mm 0mm, clip, height=14.5cm]{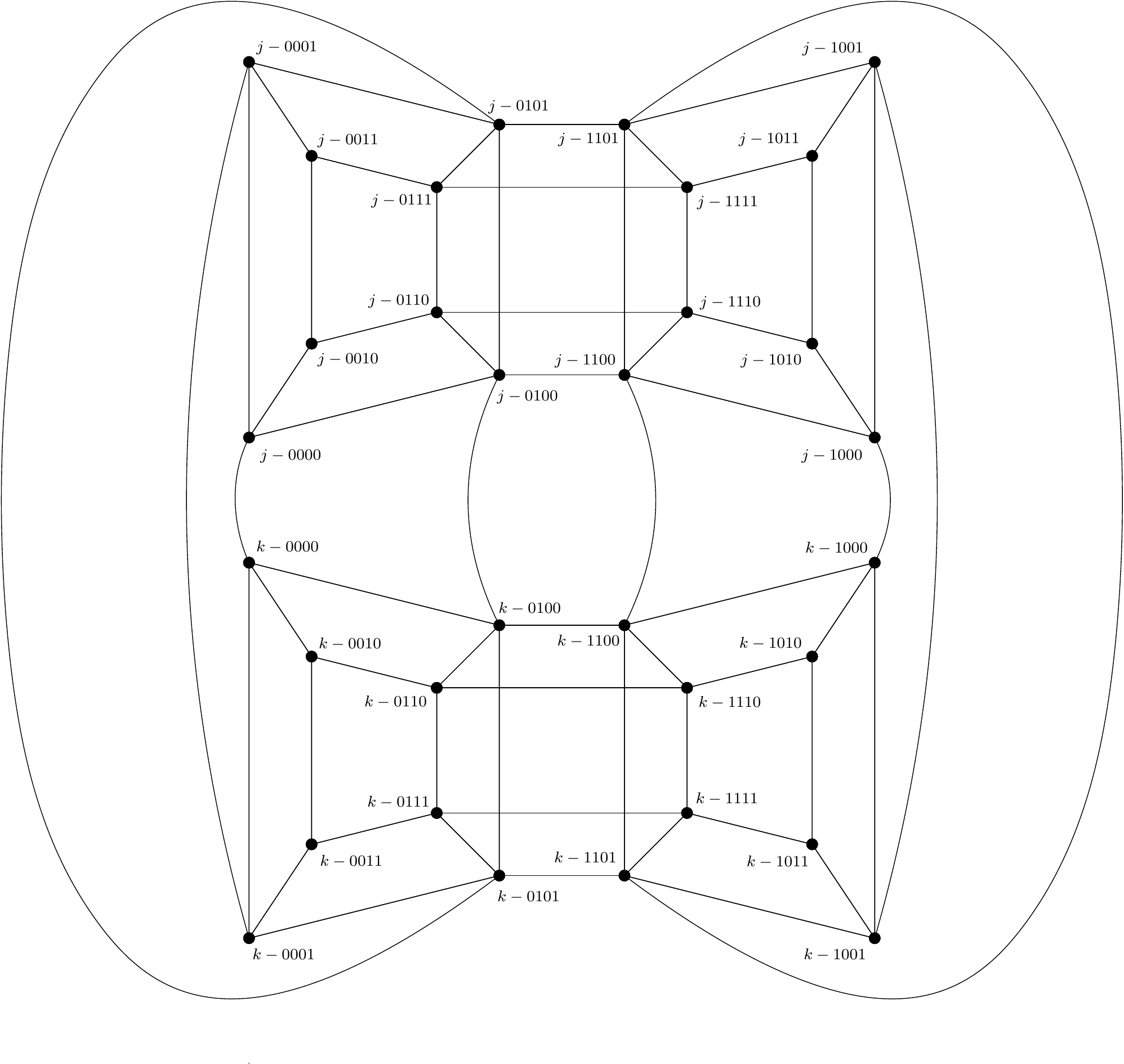} \caption{A drawing of $D_1(\hat c_1, \hat c_2)$ for $(\hat c_1, \hat c_2) \in P_1$ with eight crossings.}
\end{figure}

The four\textit{ depleted 5-cubes of Type 1} are vertex disjoint (from the form of pairs in $P_1$). We present an eight-crossing drawing of a \textit{depleted 5-cube of Type 1} in Figure \ref{F:Type1}, which proves the following claim. 

\begin{claim}
$cr(D_1(\hat c_1, \hat c_2)) \leq 8$.
\end{claim}
We similarly define $D_2(\check c_1, \check c_2)$, the \emph{Depleted 5-cube of Type 2}, according to Table \ref{Type2} given \[P_2 := \{(0111, 1111), (0101, 1101), (0100, 1100), (0110, 1110)\} \subset \binom{C_2}{2}.\]

\begin{table}[h!]
\begin{tabular}{|c|c|c|c|} 
\hline
\multicolumn{4}{|c|}{$E(D_2(\check c_1, \check c_2))$ for $\check c \in (\check c_1, \check c_2) \in P_2$.} \\
\hline
   %outer square in drawing type 1:
  $(\check{c}-b000, \check{c}-b001)$& $(\check{c}-b000, \check{c}-b100)$&  $(\check{c}-b100, \check{c}-b101)$&$(\check{c}-b001, \check{c}-b101)$\\ \hline
   %inner square in drawing type 1
   $(\check{c}-b010, \check{c}-b011)$ & $(\check{c}-b010, \check{c}-b110)$  & $(\check{c}-b110, \check{c}-b111)$ & $(\check{c}-b011, \check{c}-b111)$\\ \hline
   % diagonals in drawing type 1
   $(\check{c}-b000, \check{c}-b010)$ & $(\check{c}-b001, \check{c}-b011)$  &$(\check{c}-b100, \check{c}-b110)$ &$(\check{c}-b101, \check{c}-b111)$\\ \hline
   %Q_3 to Q_3 edges
    $(\check{c}-0011, \check{c}-1011)$ &$(\check{c}-0001, \check{c}-1001)$ &$(\check{c}-0000, \check{c}-1000)$ & $(\check{c}-0010, \check{c}-1010)$\\ 
   \hline
   % $(0111-\check c, 1111-\check c)$ & $(0101-\check c, 1101-\check c)$ & $(0100-\check c, 1100-\check c)$ & $(0110-\check c, 1110-\check c)$\\
   %\hline 
       $(\check c_1-0110, \check c_2-0110)$ & $(\check c_1-0111, \check c_2-0111)$ & $(\check c_1-0011, \check c_2-0011)$ & $(\check c_1-1011, \check c_2-1011)$\\ \hline
    $(\check c_1-1111, \check c_2-1111)$ & $(\check c_1-1110, \check c_2-1110)$ & $(\check c_1-1010, \check c_2-1010)$ & $(\check c_1-0010, \check c_2-0010)$\\
    \hline
\end{tabular}
\caption{Table of 64 edges of  \textit{depleted 5-cubes of Type 2.}}
\label{Type2}
\end{table}

\begin{figure}[h!]
\hspace*{-17mm}\includegraphics[trim = 0mm 0mm 0mm 0mm, clip, height=14.5cm]{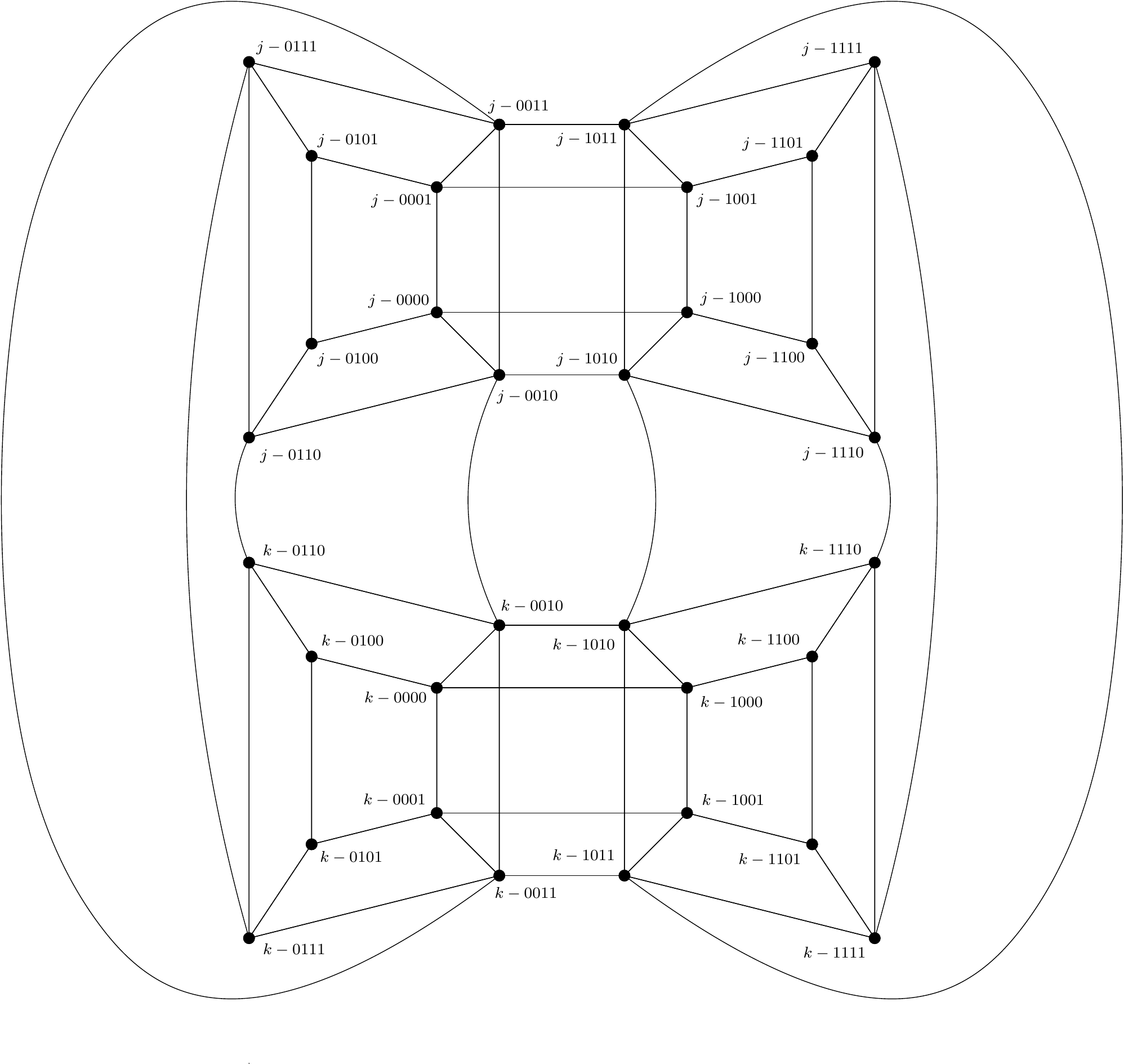}
\caption{A drawing of $D_2(\check c_1, \check c_2)$ for $(\check c_1, \check c_2) \in P_2$ with eight crossings.}
\label{F:Type2}
\end{figure}

Again, the four \textit{depleted 5-cubes of Type 2} are vertex disjoint.
An eight-crossing drawing of a \textit{depleted 5-cube of Type 2} is given in Figure \ref{F:Type2}, which proves the following claim.

\begin{claim}
$cr(D_2(\check c_1, \check c_2)) \leq 8$.
\end{claim}

Each \textit{depleted 5-cube} has 64 edges, so Plane 1 contains 512 edges. Further, no 
\textit{depleted 5-cube of Type 1} shares a vertex with a \textit{depleted 5-cube of Type 2}. This follows from the form of the pairs in $P_1$ and $P_2$ and the form of the edge sets described in Tables 1 and 2. Thus, these 512 edges can be drawn in Plane 1 with at most 64 crossings.

\begin{remark} 
Plane 2 contains all the edges of $Q_8$ which are not in Plane 1.  Moreover, $G_1 \cong G_2$.
\end{remark} 

We now provide a more illuminating description of the edges of Plane 2. The edges in Plane 2 have a symmetric representation in terms of the edges in Plane 1.  Let $\rho:E(Q_8) \to E(Q_8)$ such that \[\rho((v_pv_s, u_pu_s)) = (v_sv_p, u_su_p)\] where $v_p$ is a prefix string of length four, $v_1v_2v_3v_4$, and $v_s$ is a suffix string of length four, $v_5v_6v_7v_8$ that together define vertex $v = v_1v_2\dots v_8$.  Indeed $\rho$ captures the symmetric relationship between edges in Plane 1 and the edges in Plane 2.  Assuming an ordering on the vertices of $Q_8$ one can check that $\rho$ is indeed a bijection.  As an example, in Table \ref{Type1} we assign edge ($\hat{c}$b-000, $\hat{c}$b-001) to Plane 1.  So we send \[\rho((\hat{c}b-000,\hat{c}b-001)) = (b000-\hat{c}, b001-\hat{c})\] to Plane 2.  If we let ${\cal P}_i$ be the set of edges partitioned into Plane $i$ then ${\cal P}_2 = \rho({\cal P}_1).$ 
Moreover, the drawings provided in Figures 1 and 2 for \textit{depleted 5-cubes of Type $1$} (or \textit{Type 2}, resp.) are also drawings of their images under $\rho$. It follows that, for the edge partition we describe, each plane can be drawn with at most $64$ crossings implying that $cr_2(Q_8) \leq 128$ as desired.

A natural next step in this research is to determine whether or not this bound is sharp.  The authors believe this to be the case; however, such a proof remains elusive.  Alas, we leave the reader with the following conjecture.  

\begin{conjecture}
$cr_2(Q_8) = 128.$
\end{conjecture}

\section{Lower Bounds on \textit{structurally-symmetric} $k$-planar crossing numbers for Hypercubes}

Notably, our bi-planar drawing of $Q_8$ satisfies $G_1 \cong G_2$. This is a rather special property and is termed \emph{self-complementary} in \cite{Cza}.  It could be the case that there exists a non-isomorphic partition of $E(Q_8)$ which admits strictly fewer crossings. Yet, we wonder whether demanding that the $G_i$ be isomorphic truly forces a suboptimal number of crossings for $k$-planar drawings.  In particular, such symmetry would be expected when considering highly symmetric graphs like hyper-cubes.  

To formalize this question, we introduce the following generalization of self-complementary edge partitions.

\begin{definition} For a finite graph $G=(V,E)$, let $P$ denote an edge-partition $E = (E_1, E_2,..., E_k)$ and define $G_i=(V, E_i)$ for all $i$. If for all pairs $(r,s) \in [k] \times [k]$ we have $G_r \cong G_s$, then $P$ is a $k$\textit{-structurally-symmetric partition} of $G$.
\end{definition}

Trivially, when $|E|$ is not a multiple of $k$, no $k$\textit{-structurally-symmetric partition} of $E$ exists.\\

\begin{definition} If there exists a  $k$\textit{-structurally-symmetric partition} for $G$ that can be drawn with $cr_k(G)$ crossings then we say that the graph $G$ is $k$\textit{-structurally-symmetric}.
\end{definition}

It is unclear whether graphs exist for which any $k$\textit{-structurally-symmetric partition} of $E$ forces a sub-optimal $k$-planar drawing (which requires strictly more than $cr_k(G)$ crossings).  

In particular, we leave the reader with the following question.

\begin{question}Is the $d$-dimensional hypercube $2$\textit{-structurally-symmetric}?
\end{question}

This question motivates the following definition.

\begin{definition}
Let $cr_{kss}(G)$ denote the minimum number of crossings required among all $k$-structurally symmetric partitions of $G$.  We call $cr_{kss}$ \emph{the $k$-structurally-symmetric crossing number of $G$}.
\end{definition}

Trivially, $cr_{kss}(G)\geq cr_k(G)$. So, $k$-\textit{structurally symmetric graphs} are precisely those graphs $G$ that have $cr_k(G)=cr_{kss}(G)$. We conclude by presenting the reader questions concerning $k$-structurally-symmetric crossing numbers.

\begin{question} Characterize the set of all $k$\textit{-structurally-symmetric} graphs.  To this end, what structural properties ensure that a graph is $k$-structurally-symmteric or otherwise?
\end{question}

\begin{question}
Provide a graph for which the difference between $cr_{kss}(G)$ and $cr_k(G)$ is large (or even $>0$).  Further, is there an infinite family $(G_n)_{n \geq 1}$ such that $G_{n} \subseteq G_{n+1}$ and $(cr_{kss}(G_n) - cr_{k}(G_{n}))_{n \geq 1} \uparrow \infty$?
\end{question}

\section{Acknowledgements}
This material is based upon work that started at the Mathematics Research Communities workshop ``Beyond Planarity: Crossing Numbers of Graphs", organized by the American Mathematical Society, with the support of the National Science Foundation under Grant Number DMS 1641020.

We would like to extend our thanks to the organizers of the workshop for their commitment to engendering academic growth in young career mathematicians.  We are particularly thankful for L\'aszl\'o Sz\'ekely and his exemplary mentoring which made this project possible.

\end{document}